\documentclass[11pt]{amsart}

\usepackage{amsfonts,amsthm,amsmath,enumerate,amssymb,latexsym,color,tcolorbox,tikz,mathrsfs,bm}
\usetikzlibrary{bending}
\usepackage[text={6in,9in},centering]{geometry}
\usepackage[colorlinks,
            linkcolor=blue,
            anchorcolor=red,
            citecolor=red]{hyperref}

\graphicspath{{figures/}}

\theoremstyle{plain}
\newtheorem{theorem}{Theorem}[section]
\newtheorem{lemma}[theorem]{Lemma}

\newtheorem{proposition}[theorem]{Proposition}

\newtheorem{question}[theorem]{Question}

\theoremstyle{definition}
\newtheorem{definition}[theorem]{Definition}
\newtheorem{example}[theorem]{Example}

\theoremstyle{remark}
\newtheorem{remark}[theorem]{Remark}

\numberwithin{equation}{section}

\DeclareMathOperator{\dimh}{dim_H}
\DeclareMathOperator{\dist}{dist}

\def\R{\mathbb R}

\def\i{\bm i}
\def\j{\bm j}

\def\Q{\mathbb Q}
\def\h{\mathcal H}

\begin{document}

\title[On a self-embedding problem of self-similar sets]{On a self-embedding problem of self-similar sets}

\author{Jian-Ci Xiao}
\address{Department of Mathematics, The Chinese University of Hong Kong, Shatin, Hong Kong}
\email{jcxiao@math.cuhk.edu.hk}



\begin{abstract}
    Let $K\subset\R^d$ be a self-similar set generated by an iterated function system $\{\varphi_i\}_{i=1}^m$ satisfying the strong separation condition and let $f$ be a contracting similitude with $f(K)\subset K$. We show that $f(K)$ is relative open in $K$ if all $\varphi_i$'s share a common contraction ratio and orthogonal part. We also provide a counterexample when the orthogonal parts are allowed to vary. This partially answers a question in Elekes, Keleti and M\'ath\'e [Ergodic Theory Dynam. Systems 30 (2010)].

    As a byproduct of our argument, when $d=1$ and $K$ admits two homogeneous generating iterated function systems satisfying the strong separation condition but with contraction parts of opposite signs, we show that $K$ is symmetric. This partially answers a question in Feng and Wang [Adv. Math. 222 (2009)].
\end{abstract}

\maketitle

\section{Introduction}

For $1\leq i\leq m$, let $\varphi_i(x)=r_iO_ix+a_i$ be a contracting similitude on $\R^d$, where $0<r_i<1$ is the contraction ratio, $O_i$ is an orthogonal transformation on $\R^d$ and $a_i\in\R^d$. Let $K$ be the self-similar attractor generated by the iterated function system (or IFS in short) $\Phi:=\{\varphi_1,\ldots,\varphi_m\}$, i.e., $K$ is the unique non-empty compact set such that $K=\bigcup_{i=1}^m \varphi_i(K)$. Sometimes $\Phi$ is also referred to as a generating IFS of $K$. For a detailed introduction to self-similar sets and IFSs, please refer to the classical textbook~\cite{Fal14}.

The IFS $\Phi$ is said to satisfy the \emph{open set condition} (OSC) if there is a non-empty open set $O$ such that $\bigcup_{i=1}^m \varphi_i(O)\subset O$ and $\varphi_1(O),\ldots,\varphi_m(O)$ are pairwise disjoint; and to satisfy the \emph{strong separation condition} (SSC) when $\varphi_1(K),\ldots,\varphi_m(K)$ are pairwise disjoint. It is well known that SSC implies OSC. Under these separation conditions, the self-similar attractor $K$ has been extensively studied in many aspects, especially the fractal dimensions and self-similar measures~\cite{Hut81,Sch94,Lar95,Fal97,Fal14}. Despite this, certain geometric nature of $K$ remains mysterious even with the assumption of OSC or SSC. For example, if there is a contracting similitude $f$ that embeds $K$ into itself, what can we say about $f$? It is quite reasonable to guess that such a mapping cannot be too peculiar and there are existing results that describes this phenomenon.

In~\cite{FW09}, Feng and Wang studied this problem (or more precisely the generating IFS problem) on the real line. One of their results is the so-called \emph{logarithmic commensurability theorem}: if $\varphi_i(x)=rx+a_i$ for a common $0<r<1$, $\dimh K<1$ and the OSC is assumed, then $f(K)\subset K$ implies that $\frac{\log r_f}{\log r}\in\Q$, where $r_f$ is the contraction ratio of $f$. In other words, $r_f$ should be a rational power of $r$.

Later, Elekes, Keleti and M\'ath\'e~\cite{EKM10} generalized this result to higher dimensional cases with OSC replaced by SSC. Please see Lemma~\ref{lem:ekm} for the explicit statement. Their paper also contains many interesting findings. For example, they proved that assuming the SSC, there are only finitely many similitudes $f$ for which $f(K)\subset K$ and $f(K)$ intersects at least two level-1 cells (please see Definition~\ref{def:cell}), and for any self-similar measure $\mu$, the collection of interior points of $f(K)$ with respect to $K$ has the same $\mu$-mass as $f(K)$. 

Another related work is Algom and Hochman~\cite{AH19} which studied the self-embedding problem of Bedford-McMullen carpets. They discovered that if a carpet is created by an $m\times n$ grid with $\frac{\log m}{\log n}\notin\Q$ and is not a product set, then any similitude mapping it into itself must be an isometry composed of reflections about lines parallel to the axes.

There are also research on variants of the above embedding problem including, among others, affine embeddings of self-similar sets (\cite{Alg20}) and embeddings between self-similar or self-affine sets (\cite{DWXX11,FHR14,Alg18,FX18}).

The paper is mainly devoted to a seemingly innocent question as follows.

\begin{question}[{\cite[Question 9.3]{EKM10}}]\label{que:main}
    If $\Phi$ satisfies the SSC and $f$ is a contracting similitude mapping $K$ into itself, is it necessary that $f(K)$ is relative open in $K$?
\end{question}   

As was mentioned before, Elekes, Keleti and M\'ath\'e~\cite{EKM10} showed that $f(K)$ contains many interior points. Unfortunately, this is still far from the openness of the whole set $f(K)$. If $f(K)$ is simply a level-$n$ cell (see Definition~\ref{def:cell}) for some $n\geq 1$ then it is certainly open. Nevertheless, there are indeed many scenarios in which $f(K)$ crosses two or more cells (please refer to~\cite[Theorem 6.2]{EKM10} for an example). The way in which $f(K)$ is embedded in $K$ can be very intricate and we rarely know anything about it.

Our first result is to establish the openness of $f(K)$ when all similitudes in $\Phi$ share a common contraction ratio and orthogonal transformation. 

\begin{theorem}\label{thm:main1}
    If $r_i=r$ and $O_i=O$ for all $1\leq i\leq m$, where $0<r<1$ and $O$ is a $d\times d$ orthogonal matrix, then $f(K)$ is relative open in $K$.
\end{theorem}

Somewhat surprisingly, when the orthogonal parts are allowed to vary, we are able to construct a ``counterexample'' (please see Example~\ref{exa:main}) to the above openness statement. In other words, we have the following result. 

\begin{theorem}\label{thm:coutexa}
    Removing the assumption that $O_i\equiv O$ in Theorem~\ref{thm:main1}, $f(K)$ need not to be relative open in $K$.
\end{theorem}

We remark that in Theorem~\ref{thm:main1}, one cannot hope to relax the SSC to the OSC. For example, $\{\frac{1}{2}x,\frac{1}{2}x+\frac{1}{2}\}$ is a generating IFS of the unit interval $[0,1]$ which satisfies the OSC, but $f([0,1])$ is clearly not an open subset of $[0,1]$ for any contracting similitude $f$.

There is also a forklore open question closely related to Question~\ref{que:main}: \emph{are there two self-similar IFSs on $\R^d$ which generate the same attractor, such that one satisfies the OSC but not the SSC, and the other satisfies the SSC?} Recently, Feng, Ruan and Xiong~\cite{FRX23} gave a negative answer for homogeneous cases. More precisely, they proved that if a self-similar set admits a homogeneous generating IFS (i.e., all the similitudes in that IFS share a common contraction ratio) with OSC but not SSC, then any generating self-similar IFS of it cannot satisfy the SSC. As pointed out in their paper, an affirmative answer to Question~\ref{que:main} would lead to a negative answer to the above question. Thus, if we not only require the IFS to be homogeneous but also to satisfy the ``same orthogonal parts'' assumption, then we get a stronger conclusion (Theorem~\ref{thm:main1}); but if the orthogonal parts can vary, then one cannot expect an affirmative answer of Question~\ref{que:main}. For research on generating IFSs of self-similar sets, interested readers can refer to~\cite{FW09,DL13,DL17} and reference therein.

Finally, by using the idea in the proof of Theorem~\ref{thm:main1}, we can give a partial positive answer to another interesting question as follows.

\begin{question}[{\cite[Open Question 1]{FW09}}]
    Let $\{rx+a_i\}_{i=1}^m$ and $\{-rx+b_j\}_{j=1}^{m'}$ be two self-similar IFSs on $\R$ generating the same attractor $S$. If both of them satisfy the OSC, is it necessary that $S$ is symmetric?
\end{question}

This was answered in affirmative in~\cite{FW09} under a strong assumption called the \emph{convex open set condition}, i.e., the open set in the definition of the OSC can be chosen as a non-degenerate open interval. We are able to prove the desired symmetry under the SSC.

\begin{theorem}\label{thm:main2}
    Let $S\subset\R$ be a self-similar set. Suppose there are two generating IFSs $\{rx+a_i\}_{i=1}^m$ and $\{-rx+b_j\}_{j=1}^{m'}$ of $S$ and both of them satisfy the SSC. Then $-S=S+c$ for some $c\in\R$.
\end{theorem}

We will use $|A|$ to denote the diameter of a set $A\subset\R^d$, and  
\[
    \dist(A,B) = \inf\{\|x-y\|: x\in A, y\in B\}, \quad A,B\subset\R^d
\]
to be the distance between $A$ and $B$.

The paper is organized as follows. In Section 2 we prove Theorem~\ref{thm:main1} and construct an example to establish Theorem~\ref{thm:coutexa}. In Section 3 we prove Theorem~\ref{thm:main2}.

\section{The relative openness problem}

Let $\varphi_i,K$ be as in the beginning. Write $f(x)=r_fO_fx+a_f$, where $0<r_f<1$, $O_f$ is an orthogonal matrix and $a_f\in\R^d$. Write $\Lambda=\{1,2,\ldots,m\}$ to be the alphabet and $\Lambda^n=\{\i=i_1\cdots i_n: i_j\in\Lambda,1\leq j\leq n\}$ to be the set of words of length $n$. Also, let $\Lambda^*=\bigcup_{n=1}^\infty\Lambda^n$ denote the collection of all finite words. 

\begin{definition}\label{def:cell}
    For $n\geq 1$ and $\i=i_1\cdots i_n\in\Lambda^n$, we will refer to $\varphi_{\i}(K):=\varphi_{i_1}\circ\cdots\circ\varphi_{i_n}(K)$ as a \emph{level-$n$ cell}.
\end{definition}

Let us start with a special case where the linear part of $f$ coincide with the common one of similitudes in $\Phi$. This is actually not far away from general cases.

\begin{proposition}\label{prop:specialcase}
    Assume that $r_i=r$ and $O_i=O$ for all $1\leq i\leq m$, where $0<r<1$ and $O$ is an orthogonal transformation on $\R^d$. If we not only have $f(K)\subset K$ but also $r_f=r$ and $O_f=O$, then $f(K)$ is relative open in $K$.
\end{proposition}
\begin{proof}
    For convenience, we may assume that $|K|=1$ and write
    \[
        I = \{1\leq i\leq m: \varphi_i(K)\cap f(K)\neq\varnothing\}.
    \]
    If $\#I=1$, say $I=\{i_0\}$, then $f(K)\subset\varphi_{i_0}(K)$. We claim that $f(K)=\varphi_{i_0}(K)$. Otherwise, since $f(K)$ is a closed subset of $\varphi_{i_0}(K)$, there is some $\omega\in\Lambda^*$ such that $\varphi_{i_0\omega}(K)\cap f(K)=\varnothing$. Writing $\alpha=\dimh K$, it is well known that $0<\h^\alpha(K)<\infty$ (see~\cite[Theorem 9.3]{Fal14}). Thus 
    \[
        \h^\alpha(f(K)) \leq \h^\alpha(\varphi_{i_0}(K)\setminus\varphi_{i_0\omega}(K)) < \h^\alpha(\varphi_{i_0}(K)) = r^\alpha\h^\alpha(K) = \h^\alpha(f(K)),
    \]
    which leads to a contradiction. In particular, $f(K)$ is open as a level-$1$ cell. So it remains to consider when $\#I\geq 2$.

    Let $n\geq 1$ be so large that $r^n<\delta:=\min_{i\neq j} \dist(\varphi_i(K),\varphi_j(K))$. As a result, recalling that $|K|=1$, the diameter of every level-$n$ cell is strictly less than $\delta$. For $\mathcal{J}\subset\Lambda^{n-1}$, we call $\bigcup_{\i\in\mathcal{J}}\varphi_{\i}(K)$ a \emph{chain} if the following two conditions are met:
    \begin{enumerate}
        \item for any pair of $\i,\i'\in\mathcal{J}$, we can always find a sequence $\{\j_k\}_{k=1}^p\subset\mathcal{J}$ such that $\j_1=\i$, $\j_p=\i'$, and $\dist(\varphi_{\j_k}(K),\varphi_{\j_{k+1}}(K)) \leq r^{n-1}$ for all $1\leq k\leq p-1$.
        \item for any $\j\in\Lambda^{n-1}\setminus\mathcal{J}$, $\dist(\varphi_{\j}(K),\bigcup_{\i\in\mathcal{J}}\varphi_{\i}(K))>r^{n-1}$.
    \end{enumerate}
    Let $\mathcal{E}$ denote the collection of all chains. Under these conditions, it is not hard to see that for any $E\in\mathcal{E}$, there is exactly one $i\in I$ and $E'\in\mathcal{E}$ such that $f(E)\subset\varphi_i(E')$.

    To prove the proposition, it suffices to show that for every $E\in\mathcal{E}$, there is some $E'\in\mathcal{E}$ and some $i\in I$ such that $f(E)=\varphi_i(E')$. Suppose on the contrary that there exists $E_0\in\mathcal{E}$ not obeying this claim. Pick $i_1\in I$ and $E_1\in\mathcal{E}$ so that $f(E_0)\subset\varphi_{i_1}(E_1)$. Inductively, for each $k\geq 1$, we can find the unique $i_{k+1}\in I$ and $E_{k+1}\in\mathcal{E}$ such that $f(E_k)\subset \varphi_{i_{k+1}}(E_{k+1})$. Since $f,\varphi_{i_{k+1}}$ have the same contraction ratio, $|E_k|\leq |E_{k+1}|$ and $\h^\alpha(E_k)\leq \h^\alpha(E_{k+1})$ for all $k$. Roughly speaking, the size of these chains is ``increasing''. Also, if $f(E_k)\subsetneq \varphi_{i_{k+1}}(E_{k+1})$ then $\h^\alpha(E_k)<\h^\alpha(E_{k+1})$.
    
    Since $E_0$ does not obey our claim, $f(E_0)\subsetneq \varphi_{i_1}(E_1)$. So $\h^\alpha(E_0)<\h^\alpha(E_1)$ and hence $E_0\neq E_1$. Since there are only finitely many chains in $\mathcal{E}$, we can find $1\leq s<t$ such that $E_s=E_t$. Of course we may assume that $(s,t)$ is the earliest such pair. Note that 
    \[
        \h^\alpha(E_s) \leq \h^\alpha(E_{s+1})\leq\cdots\leq \h^\alpha(E_t)=\h^\alpha(E_s).
    \]
    Combining with the observation in the end of the last paragraph, $f(E_k)=\varphi_{i_{k+1}}(E_{k+1})$ for all $s\leq k\leq t-1$. Equivalently,
    \begin{equation}\label{eq:roekroek1}
        rOE_k+a = rOE_{k+1}+a_{i_{k+1}}, \quad \forall s\leq k\leq t-1.
    \end{equation}
    Since $f(E_{s-1})\subset \varphi_{i_s}(E_s)$ and $E_s=E_t$, we have 
    \begin{align*}
        \varphi_{i_t}(E_{s-1}) = rOE_{s-1}+a_{i_t} 
        &= (rOE_{s-1} + a) +(a_{i_t}-a) \\
        &= f(E_{s-1}) + (a_{i_t}-a) \\
        &\subset \varphi_{i_s}(E_s) + (a_{i_t}-a) \\
        &= rOE_s+a_{i_s}+a_{i_t}-a \\
        &= rOE_t+a_{i_t}+a_{i_s}-a \\
        &= rOE_{t-1}+a+a_{i_s}-a = \varphi_{i_{s}}(E_{t-1}),
    \end{align*}
    where the second last equality follows from~\eqref{eq:roekroek1}. In particular, $\varphi_{i_t}(K)\cap\varphi_{i_s}(K)\neq\varnothing$. Then $i_t=i_s$ because of the SSC assumption. As a result, $E_{s-1}\subset E_{t-1}$. But by the definition of chains, we must have $E_{s-1}=E_{t-1}$. This contradicts the ``earliest appearance'' assumption on $(s,t)$.
\end{proof}

\begin{remark}
    It is noteworthy that the above proof actually reveals more than the relative openness of $f(K)$. More precisely, it shows that $f(K)$ is a finite union of level-$n$ cells of $K$, where $n$ is the smallest integer such that $r^n<\min_{i\neq j} \dist(\varphi_i(K),\varphi_j(K))$. Furthermore, the assumption ``$O_i\equiv O$'' can be relaxed as ``$O_i=O$ for all $i\in I$''.
\end{remark}

To prove Theorem~\ref{thm:main1}, we need the following two observations.

\begin{lemma}[\cite{EKM10}]\label{lem:ekm}
    Let $g$ be a similitude on $\R^d$. If $g(K)\subset K$, then there are some integer $k\geq 1$ and two words $\i,\j\in\Lambda^*$ such that $g^k\circ\varphi_{\i}=\varphi_{\j}$.
\end{lemma}

\begin{lemma}\label{lem:openforsomepower}
    If $f^q(K)$ is an open subset of $K$ for some integer $q>1$, then so is $f(K)$.
\end{lemma}
\begin{proof}
    Otherwise, there is some $x\in K$ and a sequence $\{y_n\}_{n=1}^\infty\subset K$ such that $y_n\to f(x)$ but $y_n\notin f(K)$. So $f^{q-1}(y_n)\to f^{q-1}(f(x))=f^{q}(x)$ but $f^{q-1}(y_n)\notin f^{q-1}(f(K))=f^q(K)$, which contradicts the openness of $f^q(K)$. 
\end{proof}

\begin{proof}[Proof of Theorem~\ref{thm:main1}]
    By Lemma~\ref{lem:ekm}, we have $f^k\circ\varphi_{\i}=\varphi_{\j}$ for some $k\geq 1$ and $\i,\j\in\Lambda^*$. Let $n,n'$ be the length of $\i$ and $\j$, respectively. Comparing the linear parts, we have $r_f^{k}r^nO_f^kO^n=r^{n'}O^{n'}$. So $r_f^k=r^{n'-n}$ and $O_f^k=O^{n'-n}$. Now, applying Proposition~\ref{prop:specialcase} to $\{\varphi_\omega: \omega\in\Lambda^{n'-n}\}$ and $f^k$, we see that $f^k(K)$ is an open subset of $K$. Then the theorem follows immediately from Lemma~\ref{lem:openforsomepower}.
\end{proof}

Below we present a counterexample when the orthogonal parts vary.

\begin{example}\label{exa:main}
    For $\theta\in[0,2\pi]$, let $R_\theta=\bigl( \begin{smallmatrix} \cos\theta & -\sin\theta\\ \sin\theta & \cos\theta\end{smallmatrix}\bigr)$ be the corresponding rotation matrix as usual. We consider the IFS $\{\varphi_1,\ldots,\varphi_9\}$, where 
    \begin{equation*}
        \begin{gathered}
            \varphi_1(x)=\tfrac{1}{6}x+(-\tfrac{15}{8},\tfrac{15}{8}), \, \varphi_2(x)=\tfrac{1}{6}x+(-\tfrac{5}{2},-\tfrac{5}{4}), \, \varphi_3(x)=\tfrac{1}{6}R_{3\pi/2}x+(-\tfrac{5}{4},-\tfrac{5}{4}), \\
            \varphi_4(x)=\tfrac{1}{6}R_{\pi/2}x+(-\tfrac{5}{2},-\tfrac{5}{2}), \, \varphi_5(x)=\tfrac{1}{6}R_{\pi}x+(-\tfrac{5}{4},-\tfrac{5}{2}), \varphi_6(x)=\tfrac{1}{6}x+(\tfrac{5}{4},-\tfrac{5}{4}),\\
            \varphi_7(x)=\tfrac{1}{6}R_{3\pi/2}x+(\tfrac{5}{2},-\tfrac{5}{4}), \, \varphi_8(x)=\tfrac{1}{6}R_{\pi/2}x+(\tfrac{5}{4},-\tfrac{5}{2}), \, \varphi_9(x)=\tfrac{1}{6}R_{\pi}x+(\tfrac{5}{2},-\tfrac{5}{2}).
        \end{gathered}
    \end{equation*}
    Letting $I=[-3,3]^2$, it is easy to see that $\bigcup_{i=1}^9 \varphi_i(I)\subset I$ and hence one can use $I$ as the initial invariant set to get the attractor $K$ by iteration.
    In particular, $K\subset I$. Moreover, $\varphi_1(I),\ldots,\varphi_9(I)$ are disjoint and hence so are $\varphi_1(K),\ldots,\varphi_9(K)$. That is to say, this IFS satisfies the SSC. Please see Figure~\ref{fig:counterexample} for an illustration.

    \begin{figure}[htbp]
        \centering
        \begin{tikzpicture}[scale=0.8]
            \draw[thick,dashed] (-3,-3) rectangle (3,3);
            \draw[thick] (-3,-3) rectangle (-2,-2);
            \draw[thick] (-1.75,-3) rectangle (-0.75,-2);
            \draw[thick] (-3,-1.75) rectangle (-2,-0.75);
            \draw[thick] (-1.75,-1.75) rectangle (-0.75,-0.75);
            \draw[thick] (0.75,-3) rectangle (1.75,-2);
            \draw[thick] (2,-3) rectangle (3,-2);
            \draw[thick] (0.75,-1.75) rectangle (1.75,-0.75);
            \draw[thick] (2,-1.75) rectangle (3,-0.75);
            \draw[thick] (-2.375,1.375) rectangle (-1.375,2.375);
            \node at(-1.875,1.875) {$\varphi_1$};
            \node at(-2.5,-1.25) {$\varphi_2$};
            \node at(-1.25,-1.25) {$\varphi_3$};
            \node at(-2.5,-2.5) {$\varphi_4$};
            \node at(-1.25,-2.5) {$\varphi_5$};
            \node at(1.25,-1.25) {$\varphi_6$};
            \node at(2.5,-1.25) {$\varphi_7$};
            \node at(1.25,-2.5) {$\varphi_8$};
            \node at(2.5,-2.5) {$\varphi_9$};
            \draw[-{>[flex=1]},thick,red,dashed] (-2.5,-2.1) arc (90:180:0.4cm);
            \draw[-{>[flex=1]},thick,red,dashed] (-1.25,-2.1) arc (90:270:0.4cm);
            \draw[-{>[flex=1]},thick,red,dashed] (-1.25,-0.85) arc (90:360:0.4cm);
            \draw[-{>[flex=1]},thick,red,dashed] (1.25,-2.1) arc (90:180:0.4cm);
            \draw[-{>[flex=1]},thick,red,dashed] (2.5,-2.1) arc (90:270:0.4cm);
            \draw[-{>[flex=1]},thick,red,dashed] (2.5,-0.85) arc (90:360:0.4cm);
        \end{tikzpicture}
        \hspace*{1cm}
        \includegraphics[width=4.8cm]{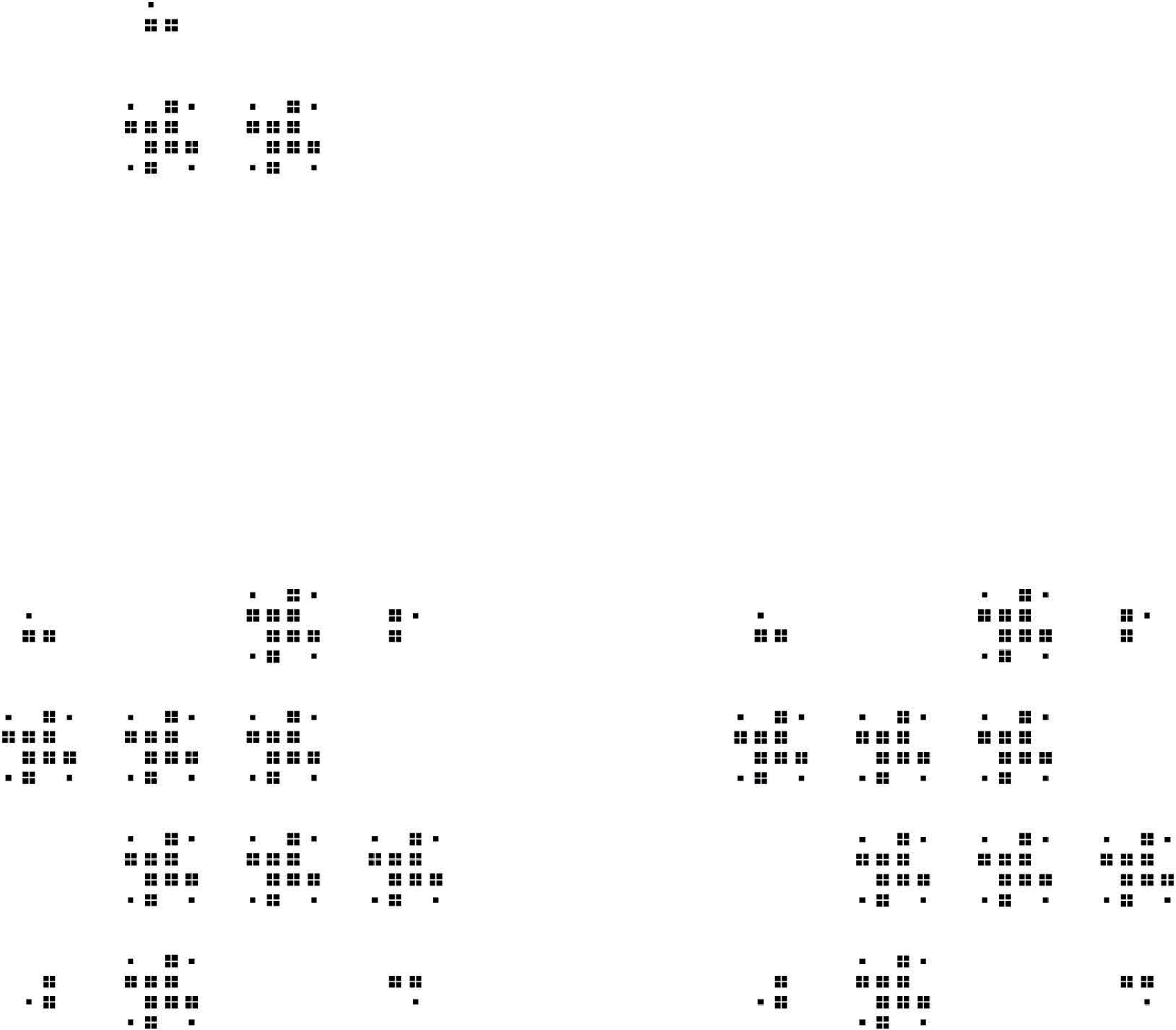}
        \caption{An illustration of the IFS and the attractor}
        \label{fig:counterexample}
    \end{figure}
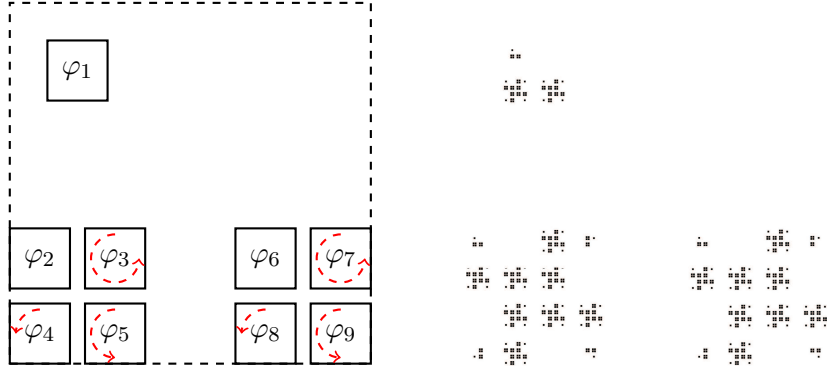

    Write $E_2=\bigcup_{i=2}^5 \varphi_i(K)$, $E_3=\bigcup_{i=6}^9 \varphi_i(K)$, $P=[-3,-\frac{3}{4}]^2$ and $Q=[\frac{3}{4},3]\times[-3, -\frac{3}{4}]$. Clearly, $E_2$ and $E_3$ only differ by a translation. A key ingredient here is the symmetry of $E_2$ (resp. $E_3$): rotating it  counterclockwise by any angle in $\{\frac{\pi}{2},\pi,\frac{3\pi}{2}\}$ at the center of P (resp. $Q$) gives us exactly the same set. Also, note that $K=\varphi_1(K)\cup E_2\cup E_3$.

    Consider $f(x)=\frac{1}{6}x+(\frac{15}{8},-\frac{15}{8})$, where $(\frac{15}{8},-\frac{15}{8})$ is nothing but the center of $Q$. We first claim that $f(K)\subset K$. A direct computation gives us the following properties:
    \begin{enumerate}
        \item $f\varphi_1=\varphi_6f$;
        \item $f(E_2)=\varphi_8(E_3)$ and $f(E_3)=\varphi_9(E_3)$ (due to the symmetry of $E_2,E_3$).
    \end{enumerate}
    By (2), writing $A=\varphi_8(E_3)\cup \varphi_9(E_3)$, we have
    \begin{align*}
        f(K) &= f(\varphi_1(K)\cup E_2 \cup E_3) \\
        &= f\varphi_1(K) \cup f(E_2) \cup f(E_3) \\
        &= f\varphi_1(K) \cup \varphi_8(E_3)\cup \varphi_9(E_3) = f\varphi_1(K) \cup A.
    \end{align*} 
    Combining with (1), 
    \begin{align*}
        f(K) = \varphi_6f(K)\cup A = \varphi_6(\varphi_6f(K)\cup A) \cup A = \varphi_6^2f(K) \cup \varphi_6(A)\cup A.
    \end{align*}
    An induction argument immediately implies 
    \begin{equation}\label{eq:fkinduction}
        f(K) = \varphi^n_6f(K) \cup \bigcup_{k=0}^{n-1}\varphi_6^k(A), \quad \forall n\geq 1.
    \end{equation}
    Note that $\dist(\varphi^n_6f(K),\bigcup_{k=0}^{n-1}\varphi_6^k(A))\leq |\varphi_6^{n-1}(K)|$, which tends to zero as $n\to\infty$. Thus
    \[
        f(K) = \overline{\bigcup_{k=0}^\infty \varphi_6^n(A)} \subset \overline{K} = K. 
    \]

    To see that $f(K)$ is not an open subset of $K$, let $x_1\in K$ be the fixed point of $\varphi_1$. By~\eqref{eq:fkinduction} and recalling the definition of $A$, we have for all $n\geq 1$ that
    \begin{align*}
        f(K) \cap \varphi_6^n\varphi_7(K) &= \Big( \varphi_6^{n+1}(K)\cup \bigcup_{k=0}^n \varphi_6^k(A) \Big) \cap \varphi_6^n\varphi_7(K) \\
        &= \big( \varphi_6^{n+1}(K) \cap \varphi_6^n\varphi_7(K) \big) \cup \bigcup_{k=0}^n \big( \varphi_6^k(A) \cap \varphi_6^n\varphi_7(K) \big) = \varnothing.
    \end{align*}
    However, since $f\varphi_1=\varphi_6f$ and $\varphi_1^n(x_1)=x_1$,
    \begin{align*}
        \dist(f(x_1),\varphi_6^n\varphi_7(K)) &\leq \dist(f\varphi_1^n(K),\varphi_6^n\varphi_7(K)) \\
        &= \dist(\varphi_6^nf(K),\varphi_6^n\varphi_7(K)) \leq |\varphi_6^n(K)|,
    \end{align*}
    which tends to $0$ as $n\to\infty$. Therefore, $f(x_1)$ is not an interior point of $f(K)$.
\end{example}

\section{The symmetry problem}

Let $S$, $\{rx+a_i\}_{i=1}^m$ and $\{-rx+b_j\}_{j=1}^{m'}$ be as in Theorem~\ref{thm:main2}. Note that the SSC implies that $m'r^s=1=mr^s$, where $s$ denotes the Hausdorff dimension of $S$. In particular, $m'=m$. To prove Theorem~\ref{thm:main2}, we need a simple observation.

\begin{lemma}\label{lem:translation}
    If there are $t_1<t_2<\cdots<t_m$ and $\alpha\in\R$ such that both of $\bigcup_{i=1}^m (S+t_i)$, $\bigcup_{i=1}^m ((-S+\alpha)+t_i)$ are disjoint unions and they are identical, then $S=-S+\alpha$.
\end{lemma}
\begin{proof}
    For notational simplicity, write $S'=-S+\alpha$. Suppose on the contrary that $S\neq S'$. Similarly as in the beginning of the proof of Proposition~\ref{prop:specialcase}, any one of $S$, $S'$ cannot be a proper subset of the other. So both of $S\setminus S'$ and $S'\setminus S$ are non-empty. Since $\bigcup_{i=1}^m (S+t_i)=\bigcup_{i=1}^m (S'+t_i)$, it is not hard to see that
    \begin{equation}\label{eq:eminusf}
        \begin{aligned}
            \bigcup_{i=1}^m ((S\setminus S')+t_i) &= \bigcup_{i=1}^m (S\setminus(S\cap S')+t_i) \\
            &= \Big( \bigcup_{i=1}^m (S+t_i) \Big) \setminus \Big( \bigcup_{i=1}^m ((S\cap S')+t_i) \Big) \\
            &= \Big( \bigcup_{i=1}^m (S'+t_i) \Big) \setminus\Big( \bigcup_{i=1}^m ((S\cap S')+t_i) \Big) = \bigcup_{i=1}^m ((S'\setminus S)+t_i),
        \end{aligned}
    \end{equation}
    where the second equality is due to the disjointness of $\{S+t_i\}_{i=1}^m$.
    
    Write $e=\inf (S\setminus S')$ and $e'=\inf (S'\setminus S)$. Then by~\eqref{eq:eminusf},
    \[
        e+t_1=\inf \Big( \bigcup_{i=1}^m ((S\setminus S')+t_i) \Big) = \inf \Big( \bigcup_{i=1}^m ((S'\setminus S)+t_i) \Big) = e'+t_1.
    \]
    So $e=e'$. Choose $\{x_k\}_{k=1}^\infty \subset S\setminus S'$ with $\lim_{k\to\infty} x_k=e$ and let $0<\delta<t_2-t_1$ be any small real number. Such a sequence exists because both of $S$, $S'$ are perfect sets. Note that when $k$ is sufficiently large, $x_k+t_1<e+t_1+\delta<e+t_2$, implying that $x_k+t_1\notin (S\setminus S')+t_j$ for all $2\leq j\leq m$. Since $e=e'$, we also have for all large $k$ that $x_k+t_1<e'+t_2$ and hence $x_k+t_1\notin (S'\setminus S)+t_j$ for all $2\leq j\leq m$. By~\eqref{eq:eminusf}, $x_k+t_1\in (S'\setminus S)+t_1$ for all large $k$. But this contradicts the fact that $x_k\in S\setminus S'$.
\end{proof}

\begin{proof}[Proof of Theorem~\ref{thm:main2}]
    Without loss of generality, we may assume the convex hull of $S$ to be $[0,1]$ and both of $\{a_i\}_{i=1}^m$, $\{b_i\}_{i=1}^m$ are increasing sequences. In particular, $a_1=0$ and $b_1=r$. For convenience, let $\varphi_i(x)=rx+a_i$ and $\psi_i(x)=-rx+b_i$. Then $a_i$, $a'_i :=b_i-r$ are the left end points of $\varphi_i(S)$, $\psi_i(S)$, respectively. Using Lemma~\ref{lem:translation}, it suffices to show that $a'_i=a_i$ for all $1\leq i\leq m$. We will prove this by induction.

    Similarly as in the proof of Proposition~\ref{prop:specialcase}, pick $n\geq 1$ so large that
    \[
        r^n<\min_{i\neq j}\dist(\varphi_i(S),\varphi_j(S)),
    \]
    adopt the definition of chains there (using this integer $n$) and denote by $\mathcal{E}$ the collection of chains of $S$. By definition, the convex hulls of these chains are disjoint. This allows us to enumerate $\mathcal{E}=\{E_1,\ldots,E_N\}$ so that $E_1,\ldots,E_N$ are located from left to right. Since every $\psi_i$ can be regarded as a self-embedding similitude of $S=\bigcup_{i=1}^m \varphi_i(S)$, similarly again as in the proof of Proposition~\ref{prop:specialcase}, for each $1\leq i\leq m$ and $1\leq t\leq N$, there are unique $i',t'$ such that $\psi_i(E_t)\subset\varphi_{i'}(E_{t'})$. 
    
    Clearly, $a'_1=b_1-r=0=a_1$. Suppose $a'_i=a_i$ for all $1\leq i\leq p$. Recall that $a_i$, $a'_i$ are left end point of $\varphi_i(S)$ and $\psi_i(S)$, respectively. Since $\psi_i(E_N)$ lies on the left end of $\psi_i(K)$, we have for all $1\leq i\leq p$ that $\psi_i(E_N)\cap\varphi_i(E_1)\neq\varnothing$. Thus $\psi_i(E_N)\subset\varphi_i(E_1)$. Also, the inclusion in turn indicates that 
    \begin{equation}\label{eq:incluandinclu}
        \begin{aligned}
            \varphi_i(E_N) = rE_N+a_i &= -(-rE_N+b_i)+a_i+b_i \\
            &= -\psi_i(E_N)+a_i+b_i \\
            &\subset -\varphi_i(E_1)+a_i+b_i = -rE_1+b_i = \psi_i(E_1).
        \end{aligned}
    \end{equation}
    In particular, $\psi_i(E_1)\cap\varphi_i(E_N)\neq\varnothing$. It then follows from the ``uniqueness'' statement in the end of the last paragraph that $\psi_i(E_1)\subset\varphi_i(E_N)$. Together with~\eqref{eq:incluandinclu}, $\varphi_i(E_N)=\psi_i(E_1)$ for $1\leq i\leq p$. As a result, $\h^s(E_1)=\h^s(E_N)$ and hence $\varphi_i(E_1)=\psi_i(E_N)$ for $1\leq i\leq p$ (due to the inclusion relationship before~\eqref{eq:incluandinclu}). 
    
    If $a'_{p+1}<a_{p+1}$, then $a'_{p+1}\subset\bigcup_{k=1}^p \varphi_i(K)$. Thus $\psi_{p+1}(E_N)\cap \bigcup_{k=1}^p \varphi_i(K)\neq\varnothing$ and hence $\psi_{p+1}(E_N)\subset \bigcup_{k=1}^p \varphi_i(K)$. So there are $1\leq k\leq p$, $1\leq s\leq N$ such that $\psi_{p+1}(E_N)\subset\varphi_k(E_s)$. A similar argument as in~\eqref{eq:incluandinclu} implies that $\psi_{p+1}(E_s)=\varphi_k(E_N)$. However, we also have by the induction hypothesis that $\psi_k(E_1)=\varphi_k(E_N)$. So 
    \begin{align*}
        \psi_k(S)\cap\psi_{p+1}(S) &\supset \psi_k(E_1) \cap \psi_{p+1}(E_s) \\
        &= \psi_k(E_1)\cap\varphi_k(E_N) = \psi_{k}(E_1) \neq\varnothing,
    \end{align*}
    which contradicts the SSC. Similarly, $a_{p+1}<a'_{p+1}$ is impossible. So $a'_{p+1}=a_{p+1}$. This completes the induction.
\end{proof}

\bigskip
\noindent{\bf Acknowledgements.}
This research is partially supported by the General Research Funds (CUHK14301017, CUHK14303021) from the Hong Kong Research Grant Council. I thank Professor De-Jun Feng for bringing the self-embedding problem to my attention and for many inspiring discussions on the subject. I am also grateful to Professor Huo-Jun Ruan for reading the manuscript and giving some comments that make this paper more readable.

\small
\bibliographystyle{amsplain}

\end{document}